\pdfoutput=1 
\documentclass[a4paper]{article}
\usepackage{amsthm,amssymb,amsmath,enumerate,graphicx}

\newtheorem{THM}{Theorem}[section]
\newtheorem{LEM}[THM]{Lemma}

\newtheorem{DEF}[THM]{Definition}

\def\shift(#1)(#2){\!\!\downarrow\!{}^{#1}_{\raise .1ex\vbox to 0pt{\vss\hbox{$\scriptstyle #2$}}}\,}
\def\ucl(#1){\lfloor #1 \rfloor}
\def\dcl(#1){\lceil #1 \rceil}
\def\specrel#1#2{\mathrel{\mathop{\kern0pt #1}\limits_{#2}}}

\def\F{\mathcal F}

\renewcommand\O{\mathcal O}

\newcommand\T{\mathcal T}

\def\lowfwd #1#2#3{{\mathop{\kern0pt #1}\limits^{\kern#2pt\raise.#3ex
\vbox to 0pt{\hbox{$\scriptscriptstyle\rightarrow$}\vss}}}}
\def\lowbkwd #1#2#3{{\mathop{\kern0pt #1}\limits^{\kern#2pt\raise.#3ex
\vbox to 0pt{\hbox{$\scriptscriptstyle\leftarrow$}\vss}}}}
\def\fwd #1#2{{\lowfwd{#1}{#2}{15}}}
\def\ve{\kern-1.5pt\lowfwd e{1.5}2\kern-1pt}
\def\vedash{{\mathop{\kern0pt e\lower.5pt\hbox{${}
     \scriptstyle'$}}\limits^{\kern0pt\raise.02ex
     \vbox to 0pt{\hbox{$\scriptscriptstyle\rightarrow$}\vss}}}}
\def\ev{\kern-1pt\lowbkwd e{0.5}2\kern-1pt}
\def\vf{\kern-2pt\lowfwd f{2.5}2\kern-1pt}
\def\vfdash{{\mathop{\kern0pt f\raise 1pt\hbox{${}
     \scriptstyle'$}}\limits^{\kern2pt\raise.02ex
     \vbox to 0pt{\hbox{$\scriptscriptstyle\rightarrow$}\vss}}}}

\def\vp{\lowfwd p{1.5}2}
\def\pv{\lowbkwd p02}
\def\vr{\lowfwd r{1.5}2}
\def\rv{\lowbkwd r02}

\def\vrdash{{\mathop{\kern0pt r\lower.5pt\hbox{${}
     \scriptstyle'$}}\limits^{\kern0pt\raise.02ex
     \vbox to 0pt{\hbox{$\scriptscriptstyle\rightarrow$}\vss}}}}
\def\rvdash{{\mathop{\kern0pt r\lower.5pt\hbox{${}
     \scriptstyle'$}}\limits^{\kern0pt\raise.02ex
     \vbox to 0pt{\hbox{$\scriptscriptstyle\leftarrow$}\vss}}}}

\def\amgis{\lowbkwd \sigma02}
\def\vs{\lowfwd s{1.5}1}
\def\sv{\lowbkwd s{1.5}1}
\def\vsdash{{\mathop{\kern0pt s\lower.5pt\hbox{${}
     \scriptstyle'$}}\limits^{\kern0pt\raise.02ex
     \vbox to 0pt{\hbox{$\scriptscriptstyle\rightarrow$}\vss}}}}
\def\svdash{{\mathop{\kern0pt s\lower.5pt\hbox{${}
     \scriptstyle'$}}\limits^{\kern0pt\raise.02ex
     \vbox to 0pt{\hbox{$\scriptscriptstyle\leftarrow$}\vss}}}}
\def\vsddash{{\mathop{\kern0pt s\lower.5pt\hbox{${}
     \scriptstyle''$}}\limits^{\kern0pt\raise.02ex
     \vbox to 0pt{\hbox{$\scriptscriptstyle\rightarrow$}\vss}}}}

\def\vsone{\lowfwd {s_1}11}
\def\svone{\lowbkwd {s_1}02}
\def\vstwo{\lowfwd {s_2}11}
\def\svtwo{\lowbkwd {s_2}02}

\def\vsidash{{\mathop{\kern0pt s_i\kern-3.5pt\lower.3pt\hbox{${}
     \scriptstyle'$}}\limits^{\kern0pt\raise.02ex
     \vbox to 0pt{\hbox{$\scriptscriptstyle\rightarrow$}\vss}}}}

\def\vS{{\hskip-1pt{\fwd S3}\hskip-1pt}} 

\def\vSr{{\vec S}_{\raise.1ex\vbox to 0pt{\vss\hbox{$\scriptstyle\ge\vr$}}}}
\def\vSstar{{\mathop{\kern0pt S\lower-1pt\hbox{$^*$}}\limits^{\kern2pt
     \vbox to 0pt{\hbox{$\scriptscriptstyle\rightarrow$}\vss}}}}
\def\vSdash{{\mathop{\kern0pt S\lower-1pt\hbox{${}
     \scriptstyle'$}}\limits^{\kern2pt\raise.1ex
     \vbox to 0pt{\hbox{$\scriptscriptstyle\rightarrow$}\vss}}}}
\def\vt{\lowfwd t{1.5}1}
\def\tv{\lowbkwd t{1.5}1}

\def\vU{{\vec U}} 

\def\es{\emptyset}
\def\sub{\subseteq}
\def\supe{\supseteq}
\def\sm{\smallsetminus}
\def\td{tree-decom\-po\-si\-tion}

\newcommand\COMMENT[1]{}
\def\?#1{\vadjust{\vbox to 0pt{\vss\vskip-8pt\leftline{%
     \llap{\hbox{\vbox{\pretolerance=-1
     \doublehyphendemerits=0\finalhyphendemerits=0
     \hsize16truemm\tolerance=10000\small
     \lineskip=0pt\lineskiplimit=0pt
     \rightskip=0pt plus16truemm\baselineskip8pt\noindent
     \hskip0pt        
     #1\endgraf}\hskip7truemm}}}\vss}}}

\title{Abstract separation systems}
 \author{Reinhard Diestel}
 \date{}

\begin{document}
\abovedisplayshortskip=-3pt plus3pt
\belowdisplayshortskip=6pt

\maketitle

\begin{abstract}\noindent
  Abstract separation systems provide a simple general framework in which both tree-shape and high cohesion of many combinatorial structures can be expressed, and their duality proved. Applications range from tangle-type duality and tree structure theorems in graphs, matroids or CW-complexes to, potentially, image segmentation and cluster analysis.

This paper is intended as a concise common reference for the basic definitions and facts about abstract separation systems in these and any future papers using this framework. 
   \end{abstract}

\section{Introduction}\label{sec:intro}

Formally, abstract separation systems are very simple objects: posets with an order-reversing involution. Think of the oriented separations $(A,B)$ of a graph, where $(A,B)\le (C,D)$ if $A\sub C$ and $B\supe D$, and the involution is given by $(A,B)\mapsto (B,A)$.

What makes such `separation systems' of graphs immediately interesting is that what little information they capture from the structure of a graph suffices to express, and to prove, two of the central theorems in graph minor theory: the tangle-tree theorem and the tangle duality theorem of Robertson and Seymour~\cite{GMX}. Even just for graphs this is not obvious.%
   \COMMENT{}
   But in fact, it can be done much more generally.

Using abstract separation systems, one can prove tangle-tree~\cite{ProfilesNew} and tangle duality~\cite{TangleTreeAbstract} theorems that apply to other combinatorial structures too. In each of these, `tangles' can encode bespoke cohesive substructures of that structure. The tangle-tree theorem shows how these can be separated in a tree-like way (thus, decomposing the given structure into its highly cohesive parts), and if there are now such `tangles' then the duality theorem shows that the structure has an overall tree shape that {\em clearly\/} cannot contain such cohesive substructures. This has already been explored in various papers~\cite{TreeSets, ProfileDuality, TangleTreeGraphsMatroids, MonaLisa} and is likely to have further applications, perhaps in areas quite different from these.

The purpose of this short paper is to collect in one place the most important definitions and basic facts about abstract separation systems, so that future papers based on them can use it as a common reference.  If more motivation is needed for why this may be a good idea, please refer to the introductions of the papers cited above, especially~\cite{ProfilesNew, TangleTreeAbstract, MonaLisa}. 

No deep theorems are proved in this paper, but many small facts are. Facts that answer obvious questions by short, but not always obvious, proofs. Small facts that can, taken together, help the reader build the kind of intuition that will make it easy to apply abstract separation systems in new contexts.

Any terminology used but not defined in this paper can be found in~\cite{DiestelBook16}.

\section{Separations}\label{sec:separations}

A \emph{separation of a set} $V$ is a set $\{A,B\}$ such that $A\cup B=V$.%
   \footnote{We can make further requirements here that depend on some structure on~$V$ which $\{A,B\}$ is meant to separate. If $V$ is the vertex set of a graph~$G$, for example, we usually require that $G$ has no edge between $A\sm B$ and $B\sm A$. But such restrictions will depend on the context and are not needed here; in fact, even the separations of a {\em set}~$V$ defined here is just an example of the more abstract `separations' we are about to introduce.}
   The ordered pairs $(A,B)$ and $(B,A)$ are its {\it orientations\/}. The {\em oriented separations\/} of~$V$ are the orientations of its separations. Mapping every oriented separation $(A,B)$ to its {\it inverse\/} $(B,A)$ is an involution%
   \COMMENT{}
   that reverses the partial ordering \[(A,B)\le (C,D) :\Leftrightarrow A\subseteq C \text{ and } B\supseteq D,\]
since the above is equivalent to $(D,C)\le (B,A)$. Informally, we think of $(A,B)$ as \emph{pointing towards}~$B$ and \emph{away from}~$A$.

More generally, a {\em separation system\/} $(\vS,\le\,,\!{}^*)$ is a partially ordered set $\vS$ with an order-reversing involution\,*. Its elements are called {\em oriented separations\/}. An {\em isomorphism\/} between two separation systems is a bijection between their underlying sets that respects both their partial orderings and their involutions.

When a given element of $\vS$ is denoted as~$\vs$, its {\em inverse\/}~$\vs^*$ will be denoted as~$\sv$, and vice versa. The assumption that * be {\em order-reversing\/} means that, for all $\vr,\vs\in\vS$,
\begin{equation}\label{invcomp}
\vr\le\vs\ \Leftrightarrow\ \rv\ge\sv.
\end{equation}

A {\em separation\/} is a set of the form $\{\vs,\sv\}$, and then denoted by~$s$. We call $\vs$ and~$\sv$ the {\em orientations\/} of~$s$. The set of all such sets $\{\vs,\sv\}\sub\vS$ will be denoted by~$S$. If $\vs=\sv$, we call both $\vs$ and $s$ {\em degenerate\/}.

When a separation is introduced ahead of its elements and denoted by a single letter~$s$, we shall use $\vs$ and~$\sv$ (arbitrarily)%
   \COMMENT{}
   to refer to its elements. Given a set $S'$ of separations, we write $\vSdash := \bigcup S'\sub\vS$%
   \COMMENT{}
   for the set of all the orientations of its elements. With the ordering and involution induced from~$\vS$, this is again a separation system.%
   \COMMENT{}%
   \footnote{When we refer to oriented separations using explicit notation that indicates orientation, such as $\vs$ or $(A,B)$, we sometimes leave out the word `oriented' to improve the flow of words. Thus, when we speak of a `separation $(A,B)$', this will in fact be an oriented separation.}

Separations of sets, and their orientations, are clearly an instance of this abstract setup if we identify $\{A,B\}$ with $\{(A,B),(B,A)\}$.

If there are binary operations $\vee$ and~$\wedge$ on our separation system~$\vS$ that make it into a lattice, i.e., such that $\vr\vee\vs$ is the supremum and $\vr\wedge\vs$ the infimum of $\vr$ and~$\vs$ in~$\vS$, we call $(\vS,\le\,,\!{}^*,\vee,\wedge)$ a {\em universe\/} of (oriented) separations. By~\eqref{invcomp}, it satisfies De~Mor\-gan's law:
\begin{equation}\label{deMorgan}
   (\vr\vee\vs)^* =\> \rv\wedge\sv.
\end{equation}%
   \COMMENT{}

A universe~$\vS$ of separations is {\em submodular\/} if it comes with a \emph{submodular order function}, a real function $\vs\mapsto |\vs|$ on~$\vS$%
   \COMMENT{}
   that satisfies $0\le |\vs| = |\sv|$ and 
 $$|\vr\vee\vs| + |\vr\wedge\vs|\le |\vr|+|\vs|$$
for all $\vr,\vs\in \vS$. We call $|s| := |\vs|$ the \emph{order} of $s$ and of~$\vs$. For every integer $k>0$, then,
 $$\vec S_k := \{\,\vs\in \vS : |\vs| < k\,\}$$
is a separation system. It need not be universe with $\lor$ and $\land$ induced by~$\vec S$,%
   \COMMENT{}
   because the supremum or infimum in~$\vec S$ of two separations in~$\vec S_k$ need not lie in~$\vec S_k$. However, by submodularity, for any two separations $\vr,\vs\in\vec S_k$ at least one of $\vr\lor\vs$ and $\vr\land\vs$ will lie in~$\vec S_k$. This motivates the following definition.

A~separation system~$\vec S$, not necessarily a universe, is called {\em submodular\/} if for any two separations $\vr,\vs\in\vec S$ either $\vr\lor\vs$ or $\vr\land\vs$ (or both) also lies in~$\vec S$.

The oriented separations $(A,B)$ of a set~$V$ form a submodular universe with respect to $|(A,B)|:= |A\cap B|$. Indeed, if $\vr = (A,B)$ and $\vs = (C,D)$, then $\vr\vee\vs := (A\cup C, B\cap D)$ and $\vr\wedge\vs := (A\cap C, B\cup D)$ are again oriented separations of~$V$, and are the supremum and infimum of $\vr$ and~$\vs$, respectively. Similarly, the oriented separations of a graph form a submodular universe of separations.

A separation $\vr\in\vS$ is {\em trivial in~$\vS$\/}, and $\rv$ is {\em co-trivial\/}, if there exists $s \in S$ such that $\vr < \vs$ as well as $\vr < \sv$.%
   \COMMENT{}
   We call such an $s$ a {\em witness\/} of $\vr$ and its triviality. If neither orientation of~$r$ is trivial, we call~$r$ {\em nontrivial\/}.%
   \COMMENT{}

The trivial oriented separations of a set~$V\!$, for example, are those of the form $\vr = (A,B)$ with $A\sub C\cap D$ and $B\supe C\cup D = V$ for some $s = \{C,D\}\ne r$ in the set $S$ considered.

Note that if $\vr$ is trivial in~$\vS$ then so is every $\vrdash \le \vr$. If $\vr$ is trivial, witnessed by~$s$, then $\vr < \vs < \rv$ by~\eqref{invcomp}. Hence if $\vr$ is trivial, then $\rv$ cannot be trivial. In particular, degenerate separations are nontrivial.%
   \COMMENT{}

\begin{LEM}\label{nontrivialexist}
If $S$ is finite,\COMMENT{}%
   \COMMENT{}
  then every trivial separation in~$\vS$ has a nontrivial witness. In particular, if $S$ is non-empty it has a nontrivial element.
 \end{LEM}

\begin{proof}
Any trivial $\vr\in\vS$ lies below a maximal trivial $\vrdash\in\vS$. If $s\in S$ witnesses the triviality of~$\vrdash$, it also witnesses that of~$\vr$. By the maximality of~$\vrdash$, neither orientation of~$s$ is trivial.%
   \COMMENT{}
\end{proof}

There can also be separations $\vs$ with $\vs < \sv$ that are not trivial.%
   \COMMENT{}
   But any\-thing smaller than these is again trivial: if $\vr < \vs\le\sv$, then $s$ witnesses the triviality of~$\vr$. Separations~$\vs$ such that $\vs\le\sv$, trivial or not, will be called {\em small\/}; note that, by~\eqref{invcomp}, if $\vs$ is small then so is every~$\vsdash\le\vs$.

The small separations $(A,B)$ of a set~$V\!$, for example, are those with $B=V$.

\begin{DEF}
 A~separation system is {\em regular\/} if it has no small elements.%
   \COMMENT{}
   It is {\em essential\/} if it has neither trivial elements nor degenerate elements. When $(\vS,\le,\!{}^*)$ is regular or essential, we also call $\vS$ and $S$ {\em regular or essential\/}.
\end{DEF}

\noindent
Note that all regular separation systems are essential.%
   \COMMENT{}

\medbreak

Since universes of separations contain infima and suprema, a finite non-empty universe is never regular, and not essential unless it consists of exactly two elements $\vr < \rv$.%
   \COMMENT{}
   The separation systems we shall study, therefore, will usually be properly contained in some universe of separations.

\begin{DEF}
The {\em essential core\/} of a separation system~$\vS$ is the essential separation system~$\vSdash$ obtained from~$\vS$ by deleting all its separations that are degenerate, trivial, or co-trivial in~$\vS$.%
   \COMMENT{}
\end{DEF}

The idea behind this definition is that $\vSdash$ is the `largest' essential separation system contained in~$\vS$. But note that this is not technically true: if $\vr$ is trivial in~$\vS$, the subsystem $\{\vr,\rv\}\sub\vS$ will also be essential but not contained in~$\vSdash$.

\medbreak

Small separations that are neither trivial nor degenerate directly precede their inverses in~$\le$:

\begin{LEM}\label{sandwich}
If $\vr$ is small and $\vr < \vs < \rv$ for some $\vs$, then $\vr$ is trivial.
 \end{LEM}

\begin{proof}
The second inequality is equivalent to $\vr < \sv$, by~\eqref{invcomp}.
\end{proof}

An essential%
   \COMMENT{}
   but irregular separation system can be made regular by deleting all pairs of the form $(\vs,\sv)$%
   \COMMENT{}
   from the relation~$\le$ viewed as a subset of~$\vS{}^2$:

\begin{DEF}
The {\em regularization\/} of an essential separation system $(\vS,\le\,,\!{}^*)$ is the triple $(\vS,\le',\!{}^*)$, where $\vr<'\!\vs$ if and only if $\vr<\vs$ and $r\ne s$.
\end{DEF}

\begin{LEM}\label{regularization}
$(\vS,\le',\!{}^*)$ is a regular separation system.
\end{LEM}

\begin{proof}
The only aspects of the assertion that $(\vS,\le,\!{}^*)$ is a separation system which rest on the existence of a pair in~$\le$  are the reflexivity and transitivity of~$\le$, and the requirement that $\sv\le\rv$ whenever $\vr\le\vs$. All these are passed on from~$\le$ to~$\le'$; for transitivity this is a consequence of Lemma~\ref{sandwich}.%
   \COMMENT{}
\end{proof}

Note that the assumption for Lemma~\ref{regularization} that $(\vS,\le,\!{}^*)$ must be essential is necessary. Indeed, as the triviality of any $\vr\in\vS$ is witnessed by some $s\ne r$, we cannot do away with it by simply `forgetting' that $\vr < \rv$: since we shall still have $\vr < \vs$ and $\vr < \sv$ as before, $\vr$~will remain trivial, and $\vr < \rv$ will continue to hold as a consequence, by~\eqref{invcomp} and transitivity of~$\le$. This is reflected in the proof of Lemma~\ref{regularization} when we use Lemma~\ref{sandwich}, which helps only when $\vr$ is small but not trivial.

The orientations of two separations $r,s$ can be related in four possible ways%
   \footnote{Actually, in eight ways; but by~\eqref{invcomp}, they come in equivalent pairs. The explicit list of relations stated here represents these four types, because every item involves~$\vr$, not~$\rv$.}%
   : as $\vr\le\vs$ or $\vr\ge\vs$ or $\vr \le \sv$ or $\vr\ge\sv$. If $r,s$ are distinct and nontrivial,%
   \COMMENT{}
   no more than one of these relations can hold:

\begin{LEM}\label{order}
If $r,s\in S$ are distinct, and have orientations~$\vr\le\vs$ such that neither $\vr$ nor $\sv$ is trivial in~$\vS$, then $\vr\not\ge\vs$ and $\vr\not\le\sv$ and $\vr\not\ge\sv$.
\end{LEM}

\begin{proof}
If $\vr\ge\vs$, then $\vr\le\vs\le\vr$ with equality, contradicting~$r\ne s$.

If $\vr\le\sv$ then $s$ witnesses that $\vr$ is trivial, contradicting our assumption.

If $\vr\ge\sv$, then $r$ witnesses that $\sv$ is trivial, contrary to assumption.
\end{proof}

A set $O\sub \vS$ of oriented separations is {\em antisymmetric\/} if $|O\cap \{\vs,\sv\}| \le 1$ for all $\vs\in\vS$: if $O$ does not contain the inverse of any of its nondegenerate elements.%
   \COMMENT{}

We call~$O$ \emph{consistent} if there are no distinct $r,s\in S$ with orientations $\vr < \vs$ such that $\rv,\vs\in O$.%
   \COMMENT{}
  By~\eqref{invcomp}, this condition is more symmetrical than it looks: if $\{\vsone,\vstwo\}$ is inconsistent, then this is witnessed by {\em both\/} $\svone < \vstwo$ and $\svtwo < \vsone$. We shall use this little fact freely whenever we need to prove consistency.%
   \COMMENT{}
   It is best remembered by thinking of an inconsistent pair of separations as pointing away from each other, an intuition which is more obviously symmetrical.

\section{Tree sets and stars}\label{sec:stars}

Two separations $r,s$ are {\em nested\/} if they have comparable orientations; otherwise they \emph{cross}. If they are nested and $\vs$ is given, then $r$ in fact has an orientation comparable with~$\vs$, and its other orientation will be comparable with~$\sv$, by~\eqref{invcomp}.

Two oriented separations $\vr,\vs$ are {\em nested\/} if $r$ and~$s$ are nested.%
   \footnote{Terms introduced for unoriented separations may be used informally for oriented separations too if the meaning is obvious, and vice versa.}%
   \COMMENT{}
   We say that $\vr$ {\em points towards\/}~$s$, and $\rv$ {\em points away from\/}~$s$, if $\vr\le\vs$ or $\vr\le\sv$.

In this informal terminology, two oriented separations are nested if and only if they are either comparable or point towards each other or point away from each other. And a set $O\sub\vS$ is consistent if and only if it does not contain orientations of distinct nested separations that point away from each other.

A~set of separations is {\em nested\/} if every two of its elements are nested. 

\begin{DEF}
\begin{enumerate}[\rm (i)]\itemsep=0pt
\item A {\em tree set\/} is a nested essential separation system. When $(\vS,\le,\!{}^*)$ is a tree set, we also call $\vS$ and $S$ {\em tree sets\/}.
\item The essential core of any%
   \COMMENT{}
   nested separation~$\vS$ is the {\em tree set induced by}~$\vS$.
\item An {\em isomorphism of tree sets} is an isomorphism of separation systems that happen to be tree sets.
\end{enumerate}
\end{DEF}

For example, the set of orientations $(u,v)$ of the edges $uv$ of a tree~$T$ form a regular tree set with respect to the involution $(u,v)\mapsto (v,u)$ and the \emph{natural partial ordering} on~$\vec E(T)$: the ordering in which $(x,y) < (u,v)$ if $\{x,y\}\ne\{u,v\}$%
   \COMMENT{}
   and the unique $\{x,y\}$--$\{u,v\}$ path in $T$ joins $y$ to~$u$. We call this tree set on $\vec E(T)$ the {\it edge tree set\/} of~$T$.

Note that a degenerate separation $s$ is never nested with another nontrivial separation~$r$, since one of its orientations, $\vr$~say, would satisfy $\vr < \vs=\sv$ and hence be trivial.%
   \COMMENT{}
   In particular, a nested separation system~$\vS$ has at most one degenerate element. If it does, then this is its only nontrivial element, and the tree set induced by~$\vS$ will be empty.

   \begin{figure}[htpb]
\centering
        \includegraphics{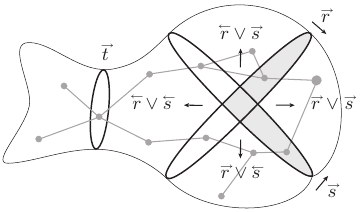}
        \caption{The separations $r$ and~$s$ cross but are nested with~$t$.}
   \label{fig:nontrans}\vskip-6pt\vskip0pt
   \end{figure}

Any two elements $\vr,\vs$ of a universe $\vS$ of separations have four \emph{corner separations},
 $\vr \lor \vs$, $\rv \lor \sv$, $\vr \lor \sv$ and~$\rv \lor \vs$, that are nested with both $\vr$ and~$\vs$. Note that the set of these four depends only on $r$ and~$s$, not on their orientations.%
   \COMMENT{}
   Figure~\ref{fig:nontrans} indicates an example for graph separations.

Figure~\ref{fig:nontrans} also illustrates the often-used {\em fish lemma\/}:

\begin{LEM} \label{lem:fish}
   Let $r,s\in S$ be two crossing separations. Every separation $t$ that is nested
   with both $r$ and $s$ is also nested with all four corner separations of $r$ and $s$.
\end{LEM}

\begin{proof}
Since $t$ is nested with $r$ and~$s$, it has an orientation pointing towards~$r$, and one pointing towards~$s$. If these orientations of $t$ are not the same, then $\vr\le\vt\le\vs$ for suitable orientations of $r,s,t$. In particular, $r$ and $s$ are nested, contrary to our assumption. Hence $t$ has an orientation $\vt$ that points towards both $r$ and~$s$.

Now $r$ and~$s$ have orientations $\vt\le\vr$ and $\vt\le\vs$. Since $\wedge$ and $\vee$ denote infima and suprema in~$\vS$, we have $\vt\le\vr\land\vs = (\rv\lor\sv)^*$ by~\eqref{deMorgan}, as well as trivially $\vt\le\vr\lor\vs$ and $\vt\le\vr\lor\sv$ and $\vt\le\rv\lor\vs$.
   \end{proof}

\begin{DEF}
A {\em star (of separations)\/} is a set $\sigma$ of nondegenerate oriented nested separations whose elements point towards each other:  $\vr\le\sv$ for all distinct $\vr,\vs\in\sigma$.
\end{DEF}

\begin{figure}[htpb]
\centering
        \includegraphics{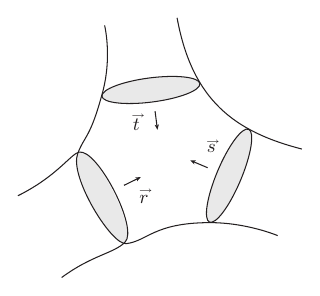}
        \caption{A star of three set separations $\vr,\vs,\vt$.}
   \label{fig:star}\vskip-6pt\vskip0pt
   \end{figure}

 We allow $\sigma=\es$. Note that stars of separations are nested. They are also consistent: if distinct $\rv,\vs$ lie in the same star we cannot have $\vr < \vs$, since also $\vs\le \vr$ by the star property. Figure~\ref{fig:star} shows a star of three set separations.

A star $\sigma$ is {\em proper\/} if, for all distinct%
   \COMMENT{}
   $\vr,\vs\in\sigma$, the relation $\vr\le \sv$ required by the definition of `star' is the only one among the four possible relations between orientations of distinct $r$ and~$s$: if $\vr \le \sv$ but $\vr\not\le\vs$ and $\vr\not\ge\vs$ and $\vr\not\ge\sv$.%
   \COMMENT{}

\begin{LEM}\label{properstars}
\begin{enumerate}[$\!\!$\rm(i)]\itemsep=0pt
\item $\!$A nested set of nondegenerate%
   \COMMENT{}
   oriented separations is a proper star if and only if it is an antisymmetric, consistent antichain.
\item $\!$A star is proper if and only if it is an antisymmetric anti\-chain.%
   \COMMENT{}
       \qed
\end{enumerate}
\end{LEM}

Let us call a star $\sigma\in\vS$ {\em proper in~$\vS$} if it is proper and is not a singleton~$\{\sv\}$ with $\sv$ co-trivial in~$\vS$. We shall call such stars {\em co-trivial singletons\/}.

Non-singleton proper stars cannot contain co-trivial separations. In fact, they cannot contain any separation whose inverse is small:

\begin{LEM}\label{cotrivialstar}
If a proper star $\sigma\sub\vS$ has an element~$\vs$ such that $\sv$ is small, then it has no other elements. In particular, if $\sigma$ is proper in~$\vS$ then none of its elements is co-trivial in~$\vS$.
\end{LEM}

\begin{proof}
Suppose there exists some $\vr\in\sigma\sm\{\vs\}$. Then $\vr\le\sv\le\vs$, since $\sv$ is small. Hence $\sigma$ is not an antichain and thus not a proper star.
\end{proof}

The simplest example of an improper star is one violating antisymmetry, e.g., the star $\{\vs,\sv\}$ for nondegenerate~$s$. Such a star may contain further separations~$\vr$, but note that these must be trivial, as witnessed by~$s$. If a star is antisymmetric but fails to be an antichain, containing separations $\vr < \vs$ say, then again $\vr$ must be trivial, since we also have $\vr\le\sv$ by the star property (and $\vr\ne\sv$ by antisymmetry, so $r\ne s$).%
   \COMMENT{}
   By Lemma~\ref{properstars}\,(ii), we thus have

\begin{LEM}\label{2stars}
Any improper star in $\vS$ that is not a co-trivial singleton and is not of the form~$\{\vs,\sv\}$ contains a separation that is trivial in~$\vS$.%
   \COMMENT{}
   \qed
\end{LEM}

Our partial ordering on~$\vS$ also relates its subsets, and in particular its stars: for $\sigma,\tau\sub\vS$ we write $\sigma\le\tau$ if for every $\vs\in\sigma$ there exists some $\vt\in\tau$ with $\vs\le\vt$.%
   \COMMENT{}
    This relation is obviously reflexive and transitive, but in general it is not antisymmetric: if $\sigma$ contains separations $\vs < \vt$, then for $\tau = \sigma\sm\{\vs\}$ we have $\sigma < \tau < \sigma$ (where $<$ denotes `$\le$ but not~=').%
   \COMMENT{}
   However, it is antisymmetric on antichains, and thus in particular on proper stars:

\begin{LEM}\label{antichains}
The set of all proper stars $\sigma\sub\vS$%
   \COMMENT{}
   is partially ordered by~$\le$.
\end{LEM}

\begin{proof}
We only have to show antisymmetry. If $\sigma\le\tau\le\sigma$, then for every $\vs\in\sigma$ there are $\vt\in\tau$ and $\vsdash\in\sigma$ such that $\vs\le\vt\le\vsdash$. If $\sigma$ is an antichain this implies $\vs=\vsdash$, and hence $\sigma\sub\tau$. Likewise $\tau\sub\sigma$, so $\sigma = \tau$.
\end{proof}

When we speak of {\em maximal proper stars in\/} a separation system~$(\vS,\le\,,\!{}^*)$, we shall always mean stars that are $\le$-maximal in the set of stars that are proper in~$\vS$. These stars $\sigma$ need not be maximal among all the stars in~$\vS$, not even among the proper ones;%
   \COMMENT{}
   for example, there may be a co-trivial singleton~$\{\rv\}$ in~$\vS$ (which is a proper star, though not proper in~$\vS$) such that $\vs\le\rv$ for all $\vs\in\sigma$.%
   \COMMENT{}\looseness=-1

Lemma~\ref{cotrivialstar} tells us that proper stars in a separation system~$\vS$ cannot contain separations that are co-trivial in~$\vS$. Our next lemma says that maximal proper stars in nested%
   \COMMENT{}
   separation systems~$\vS$ do not contain trivial separations either, and will thus lie in the tree set that $\vS$ induces:

\begin{LEM}\label{maximalproperstars}
Maximal proper stars in a nested separation system~$\vS$ without degenerate elements%
   \COMMENT{}
   contain no separations that are trivial or co-trivial in~$\vS$.%
   \COMMENT{}
   \end{LEM}

\begin{proof}
Let $\sigma$ be a proper star in~$\vS$. By Lemma~\ref{cotrivialstar} it contains no co-trivial separations. We assume that $\sigma$ has an element~$\vr$ that is trivial in~$\vS$, witnessed by $r'\in S$ say, and show that $\sigma$ is not maximal among the proper stars in~$\vS$.

Not both orientations of~$r'$ can be trivial;%
   \COMMENT{}
   let $\rvdash$ be one that is not. Let $\sigma' > \sigma$ be obtained from $\sigma$ by replacing $\vr$ with~$\vrdash$ and deleting any $\vs\in\sigma$ such that $\vs < \vrdash$. Let us show that $\sigma'$~is again a proper star in~$\vS$; this will show that $\sigma$ was not maximal among these.

No $\vs\in\sigma'\sm\{\vrdash\}$ can point away from~$r'$: that would imply either $\vr < \vrdash\le\vs$ or $\vr < \rvdash\le\vs$, contradicting the fact that $\vr$ and~$\vs$ lie in the same proper star~$\sigma$. Hence every such $\vs$ points towards~$r'$, because $S$ is nested. But $\vs\not < \vrdash$ by definition of~$\sigma'$,%
   \COMMENT{}
   so $\vs\le\rvdash$. Thus, $\sigma'$~is indeed a proper star.%
   \COMMENT{}
   It is also proper in~$\vS$, because $\vrdash$ was chosen not to be co-trivial in~$\vS$.%
   \COMMENT{}
   \end{proof}

\section{Orientations of separation systems}\label{sec:orientations}

An \emph{orientation} of a separation system~$\vS$,%
   \COMMENT{}
   or of a set~${S}$ of separations, is a set $O\sub{\vS}$ that contains for every $s\in{S}$ exactly one of its orientations $\vs,\sv$. A~\emph{partial orientation} of~${S}$ is an orientation of a subset of~${S}$: an antisymmetric subset of~$\vS$.

We shall be interested particularly in consistent orientations of separation systems~$\vS$. If $\vS$ comes from a concrete combinatorial structure which its elem\-ents separate, its consistent orientations can be thought of as pointing to the locations in this combinatorial structure which the separations in~$S$ separate from each other. Examples include all the classical highly connected substructures of graphs and matroids that have been studied in the context of width parameters, such as blocks, tangles or brambles~\cite{duality1inf, confing, ProfileDuality, TangleTreeGraphsMatroids, ProfilesNew}. They also include the vertices and the ends of a tree: these correspond to the consistent orientations of its edge set~$E$, since these always point towards a unique vertex or end~\cite{TreeSets}.\looseness=-1

 While these locations might originally be identified by concrete substructures of a structure~$X$ which our separation system~$\vS$ separates, we now see that $X$ is not in fact needed to identify them: as the locations are identified by consistent orientations of~$\vS$, and consistency is defined in terms of~$\vS$ alone, without reference to~$X$, these locations too can be described purely in terms of~$\vS$, with reference only to the axiomatic properties of~$\vS$.%
   \COMMENT{}

This shift of paradigm enables us to treat all kinds of locations in discrete structures uniformly. For example, we can prove very general duality theorems asserting that a separation system $\vS$ either admits certain consistent orientations (which for concrete choices of $\vS$ might correspond to certain types of highly cohesive substructures of a given structure) or else contains a tree set whose consistent orientations all point to locations that are `small' in a corresponding dual sense~\cite{TangleTreeAbstract}, and clearly too small to accommodate such a highly cohesive substructure.

But our general framework also defines new kinds of `locations' in discrete structures that had not previously been studied: every consistent orientation of a natural separation system of a given structure can in principle be thought of as such a `location'.

To emphasise this point, and to support our intuition, we shall therefore think of the consistent orientations also of abstract separation systems~$\vS$ in this way: as regions of some unknown combinatorial structure which $S$ `separates', regions that are either too small or too highly connected to be split by the separations in~$S$.

\medbreak

Every consistent orientation $O$ of a separation system~$\vS$ is {\em closed down\/} in~$\vS$: if $\vr,\vs\in\vS$ satisfy $\vs\in O$ and $\vr\le\vs$, we also have $\vr\in O$, since otherwise $\rv\in O$, contradicting the consistency of~$O$. Conversely, every orientation of $\vS$ that is closed down in the ordering of~$\vS$ is obviously consistent.

A consistent orientation of a separation system~$\vS$ cannot contain any separations that are co-trivial in~$\vS$: if $\vr$ is trivial, witnessed by $s\in S$, say, then $s$ cannot be oriented consistently with~$\rv$. The following lemma provides a kind of converse to this observation:

\begin{LEM}\label{treenodes} Let $P$ be a consistent partial orientation of a separation system~$\vS$.
\begin{enumerate}[\rm(i)]\itemsep=0pt
   \item $P$ extends to a consistent orientation $O$ of~$S$ if and only if no element of~$P$ is co-trivial in~$\vS$.
   \item Given any maximal element $\vp$ of~$P$, the orientation $O$ in~{\rm(i)} can be chosen with $\vp$~maximal in~$O$\penalty-200\ if and only if $p$ is nontrivial in~$S$.%
   \COMMENT{}
   \item If $S$ is nested, then the orientation $O$ in {\rm(ii)} is unique.
\end{enumerate}
\end{LEM}

\begin{proof}
(i) The forward implication follows from the fact that no consistent orientation of $S$ can contain a co-trivial separation~$\rv$: if the triviality of $\vr$ is witnessed by~$s\in S$, then both $\{\rv,\vs\}$ and $\{\rv,\sv\}$ would be inconsistent.

For the backward implication, use Zorn's lemma to extend $P$ to a maximal consistent partial orientation $O$ of~$S$. Suppose $S$ has an element $s$ such that neither $\vs$ nor $\sv$ lies in~$O$. Then $O\cup\{\vs\}$ is inconsistent, which means that there is some $\rv\in O$ such that $\vr < \vs$.%
   \COMMENT{}
   And $O\cup\{\sv\}$ is inconsistent, so there is some $\vt\in O$ such that $\vs < \vt$. But then $\vr < \vt$ with $\rv,\vt\in O$, contradicting the consistency of~$O$.

(ii) The forward implication is again clear. Indeed, if $\vp < \vs$ and $\vp < \sv$, then $\vp$ cannot be maximal in any orientation of a set containing~$s$. Since, by~(i), $\vp$~cannot be co-trivial, this proves that $p$ is nontrivial.

For the backward implication let $P'$ be the set of separations $\vs\in\vS$ such that $\vs$ points towards~$p$ but $\sv$ does not. Neither $\vp$ nor $\pv$ lies in~$P'$.%
   \COMMENT{}
   Let us show that $P\cup P'$ is still consistent. If $\sv,\vsdash\in P'$ are inconsistent, then they both point towards~$p$ and $\vs < \vsdash$. But then $\vs$ also points towards~$p$, contradicting the definition of~$P'$. If $\rv\in P$ is inconsistent with $\vs\in P'$, then $\vr < \vs\le\vp$ or $\vr < \vs\le\pv$.%
   \COMMENT{}
   The first of these cases contradicts the consistency of~$P$ (which contains both $\rv$ and~$\vp$).%
   \COMMENT{}
   In the latter case we have $\vp \le \sv < \rv\in P$, contradicting the maximality of $\vp$ in~$P$. This completes the proof that $P\cup P'$ is consistent.

By assumption~(i), no element of $P$ is co-trivial in~$\vS$. Let us show that also no $\vs\in P'$ is co-trivial in~$\vS$. If it is, then $\sv$ is trivial, and hence $\vp\not\le\sv$ since $\vp$ is non-trivial by assumption. The fact that $\vs$ points to~$p$ thus means that $\vs < \vp$ (since we cannot have $\vs \le \pv$). But then our assumed co-triviality of $\vs$ implies that $\vp$ too is co-trivial. This contradicts the fact that $\vp\in P$.

By~(i), we can extend $P\cup P'$ to a consistent orientation $O$ of~$S$. To show that $\vp$ is maximal in~$O$, assume there exists $\vs\in O$ such that $\vp < \vs$. Then $\sv$ points towards~$p$.%
   \COMMENT{}
    But $\vs$ does not: we cannot have $\vs\le\vp$ (as $\vp < \vs$), we cannot have $\vs = \pv$ (since $\vp\in O$ implies $\pv\notin O$), and we cannot have $\vs < \pv$, since then $\vp < \sv$ as well as $\vp < \vs$, contradicting the non-triviality of~$\vp$. Hence $\sv\in P'$. But then $O\supe P'$ contains $\sv$ as well as~$\vs$, a contradiction.

(iii) Let $\vp$ be as in~(ii), and consider any $\vs\in\vS$ such that $\vs < \vp$. Then $\{\sv,\vp\}$ is inconsistent, so $\sv\notin O$ and hence $\vs\in O$. Now consider any $\vs\in\vS$ such that $\vp < \vs$. Then $\vs\in O$ would contradict the maximality of $\vp$ in~$O$, so $\sv\in O$. Hence no matter how the $O$ in (ii) was chosen, it orients every $s\in S$ that is nested with $p$ in the same way.
   \end{proof}

In order to identify a consistent orientation of~$\vS$, we need to know only its maximal elements. Indeed, given a subset $\sigma\sub\vS$, let us write
 $$\dcl(\sigma) := \{\, \vr\in\vS\mid\exists \vs\in\sigma\colon \vr\le\vs\,\}$$
 for its {\em down-closure\/} in~$\vS$. If $\sigma$ is antisymmetric, we let
 $$\amgis:=\{\sv\mid\vs\in\sigma\}\sm\sigma.$$%
   \COMMENT{}
   Note that $\amgis = \{\sv\mid\vs\in\sigma\}$ if $\sigma$ has no degenerate element, e.g., if~$\sigma$ is a star.%
   \COMMENT{}

If $S$ is finite, we can recover any consistent orientation of~$S$ from the set of its maximal elements. More specifically:

\begin{LEM}\label{maxelts}
Let $\sigma$ be the set of maximal elements of a consistent orientation~$O$ of a separation system~$\vS$%
   \COMMENT{}
   such that every element of $O$ lies below some element of~$\sigma$. Then $O = \dcl(\sigma)\sm\amgis.$ In~particular, $O$~is uniquely determined by~$\sigma$.

If $\vS$ is regular, then $O = \dcl(\sigma)$.
\end{LEM}

\begin{proof}
Every $\vr\in O$ lies in~$\dcl(\sigma)\sm\amgis$,%
   \COMMENT{}
   by definition of~$\sigma$. Conversely if $\vr\in \dcl(\sigma)\sm\amgis$, then either $\vr\in\sigma\sub O$ or $\vr < \vs\in\sigma\sub O$ for some $\vs\ne\rv$. In the latter case we have $\rv\notin O$ since $O$ is consistent, and hence again $\vr\in O$.

If $\vS$ is regular, then no element of $\amgis = \{\sv\mid\vs\in\sigma\}$ lies in~$\dcl(\sigma)$.%
   \COMMENT{}
   For if $\vs\in\sigma$ and $\sv\le\vt\in\sigma$ then $\sv<\vt$ since $\{\vs,\sv\}\not\sub O$, which makes $\{\vs,\vt\}\sub O$ inconsistent (contradiction) unless $s=t$. But then $\sv<\vt=\vs$, contra\-dicting the regularity of~$\vS$.
   \end{proof}

\begin{LEM}\label{ConsTreeSets}
The consistent orientations of a separation system $\vS$ are precisely the consistent orientations of its essential core~$\vSdash$ together with all its trivial separations and the unique orientation of any degenerate separation in~$S$.
  \end{LEM}

\begin{proof}
Let $\vec R$ be the set of trivial or degenerate separations in~$\vS$. Every consistent orientation $O$ of~$\vS$ contains~$\vec R$, by Lemma~\ref{treenodes}(i), and $O\cap\vSdash$ is a consistent orientation of~$\vSdash$.

Conversely, if $O'$ is a consistent orientation of~$\vSdash$ then $O'\cup\vec R$ is a consistent orientation of~$\vS$. Indeed, no two elements of $\vec R$ can point away from each other, since this would make both their inverses trivial;%
   \COMMENT{}
   but neither a trivial nor a degenerate separation has a trivial inverse.%
   \COMMENT{}
   But neither can $\vr\in\vec R$ and $\vsdash\in O'$ point away from each other, since this would make $\svdash < \vr$ trivial,%
   \COMMENT{}
   contradicting the fact that $s'\in S'$ is nontrivial. Hence no two separations in $O'\cup\vec R$ point away from each other, which shows that $O'\cup\vec R$ is consistent.
\end{proof}

The consistent orientations of tree sets will be of particular interest to us.
Let us say that a subset $\sigma$ of a nested separation system~$\vS$ {\em splits\/}~$\vS$ if $S$ has a consistent orientation $O$ whose set of maximal elements is precisely~$\sigma$ and which satisfies $O\sub\dcl(\sigma)$.%
   \COMMENT{}%
   \COMMENT{}
   For example, the edge tree set of a tree is split precisely by the sets~$\vec F_t$ of edges at a node~$t$, oriented towards~$t$ \cite{TreeSets}.

The term `split' comes from the fact that only trivial separations $\vr$ can lie below more than one element of a splitting set~$\sigma$. Indeed, if $\vs$ and~$\vsdash$ are both maximal in~$O$ and consistently nested with each other, they must point towards each other. Hence $\vr\le\vs < \svdash$%
   \COMMENT{}
   as well as $\vr\le\vsdash < \sv$. But this makes $\vr$ trivial unless $\vr=\vs$ or $\vr=\vsdash$ (and $\vr$ is small).%
   \COMMENT{}

Splitting subsets contain no trivial separations:

\begin{LEM}\label{SplitTreeSet}
Let $\vS$ be a nested separation system.
\begin{enumerate}[\rm(i)]\itemsep=0pt
 \item If $S$ has a degenerate element~$s$, then $\{\vs\}$ is the unique splitting set in~$\vS$.
 \item If $S$ has no degenerate element, then the subsets splitting~$\vS$ are precisely the subsets that split the tree set which $\vS$ induces.
\end{enumerate}
\end{LEM}

\begin{proof}
(i) If $S$ has a degenerate element~$s$, then $s$ is its only nontrivial element,%
   \COMMENT{}
   so the tree set which $\vS$ induces is empty. By Lemma~\ref{ConsTreeSets}, therefore, $S$ has a unique consistent orientation~$O$, which consists of all its trivial elements and~$\vs=\sv$. Then $O=\dcl(\vs)$ by Lemma~\ref{nontrivialexist}.%
   \COMMENT{}
   So $\{\vs\}$ splits~$\vS$, but no other subset of~$\vS$ does.

(ii) Consider any element $\vp$ of a splitting subset $\sigma$ of~$\vS$. Since $\vp$ is maximal in a consistent orientation of~$S$, Lemma~\ref{treenodes}(ii) applied with $P=\{\vp\}$ implies that $\vp$ is not trivial in~$\vS$. Hence $\sigma$ lies in the tree set $\vSdash$ that $\vS$ induces,%
   \COMMENT{}
   which it clearly also splits.%
   \COMMENT{}

Conversely, if $\sigma'$ splits~$\vSdash$, witnessed by the orientation $O'$ of~$\vSdash$, say, then, by Lemma~\ref{nontrivialexist},%
   \COMMENT{}
   adding to $O'$ the trivial separations of~$\vS$ extends it to a consistent orientation of~$\vS$ whose set of maximal elements is still~$\sigma'$,%
   \COMMENT{}
   and which splits~$\vS$.
\end{proof}

In our example where $\vS$ is the edge tree set of a tree~$T$, the subsets $\vec F_t$ of~$\vS$ that split it can be described just in terms of $\vS$ itself, without reference to orientations (unlike in the definition of `split'): they are simply the maximal proper stars in~$\vS$~\cite{TreeSets}.%
   \COMMENT{}

For arbitrary tree sets~$\vS$, this remains true with one curious exception. But this exception, too, can be described without reference to orientations of~$\vS$. Consider the partial ordering $\le$ of the proper stars in~$\vS$ from Lemma~\ref{antichains}.

\begin{LEM}\label{Ostars}
Let $\vS$ be a nested separation system%
   \COMMENT{}
    without a degenerate element.%
   \COMMENT{}
 A~subset $\sigma\sub\vS$ splits~$\vS$ if and only if either
    \begin{enumerate}[\rm (i)]\itemsep=0pt\vskip-3pt\vskip0pt
\item $\sigma$ is a maximal proper star in~$\vS$;%
   \COMMENT{}
    or
\item $\sigma$ is a proper star in~$\vS$ that contains a small separation, and every $\sigma'>\sigma$ that is a proper star in~$\vS$%
   \COMMENT{}
   is of the form $\sigma'=\{\sv\}$ with $\vs\in\sigma$ small.%
   \COMMENT{}
  \end{enumerate}
\end{LEM}

\begin{proof}
We assume first that $\sigma$ splits~$\vS$, and show that it satisfies (i) or~(ii). Let $\sigma$ be the set of maximal elements of the consistent orientation~$O$ of~$S$.

By Lemma~\ref{properstars}\,(i), $\sigma$ is a proper star.%
   \COMMENT{}
  It is also a proper star in~$\vS$, since $O$ is consistent and so $\sigma\sub O$ cannot be a co-trivial singleton in~$\vS$ (cf.\ Lemma~\ref{treenodes}\,(i)). Suppose $\sigma' > \sigma$ is another proper star in~$\vS$. We shall prove that unless we obtain a contradiction, which will establish~(i), we have (ii) witnessed by~$\sigma'$.

Since $\sigma'\not\le\sigma$ (cf.\ Lemma~\ref{antichains}),%
   \COMMENT{}
   there exists $\vsdash\in\sigma'$ such that $\vsdash\not\le\vs$ for all~$\vs\in\sigma$. Since every element of~$O$ lies below some element of~$\sigma$,%
   \COMMENT{}
   this means that $\vsdash\notin O$, and hence $\svdash\in O$. Hence there exists $\vs\in\sigma$, and as $\sigma\le\sigma'$ also some $\vsddash\in\sigma'$, such that $\svdash \le \vs\le\vsddash$. As $\vsdash,\vsddash\in\sigma'$, this contradicts the fact that $\sigma'$ is a proper star, establishing~(i)~-- as long as $\vsdash\ne\vsddash$.

If $\vsdash=\vsddash$, we have $\svdash\le\vs\le\vsdash$.%
   \COMMENT{}
   Applying \eqref{invcomp} to the second of these inequalities we obtain $\svdash\le\sv$ (as well as $\svdash\le\vs$, the first inequality). If $s\ne s'$ this means that $\svdash$ is trivial, contradicting Lemma~\ref{cotrivialstar} for~$\sigma'$.%
   \COMMENT{}
   Hence $s=s'$. Since $\vsdash\ne\vs$ by the choice of~$\vsdash$, this means that $\vsdash = \sv$. Our double inequality now yields $\vs\le\sv$, so $\vs$ is small, and $\sigma'$ has no other element than $\vsdash=\sv$ by Lemma~\ref{cotrivialstar}.

Conversely, assume that $\sigma$ satisfies (i) or~(ii). Like all stars, $\sigma$~is consistent. By Lemmas \ref{cotrivialstar} and~\ref{treenodes}\,(i) we can extend~$\sigma$ to a consistent orientation $O$ of~$S$. We shall prove that $O\sub\dcl(\sigma)$: then $\sigma$ clearly splits~$\vS$.
   \COMMENT{}

Suppose there exists $\vsdash\in O\sm\dcl(\sigma)$. Let $\sigma':= (\sigma\cup\{\vsdash\})\sm (\dcl(\vsdash)\cap\sigma)$. Since $\sigma'$ is contained in the consistent and antisymmetric set~$O$, and is an antichain because $\sigma$ is proper, Lemma~\ref{properstars} implies that $\sigma'$ is a proper star.%
   \COMMENT{}
   It is even proper in~$\vS$, since $\vsdash\in O$ cannot be co-trivial in~$\vS$ because $O$ is consistent.

As $\sigma < \sigma'$ by definition of~$\sigma'$,%
   \COMMENT{}
   we thus cannot have~(i). So we must have~(ii), with $\sigma'= \{\sv\}$ for some $\vs\in\sigma$. As both $\vs\in\sigma$ and $\sv\in\sigma'$ lie in~$O$, which is antisymmetric, this means that $s$ is degenerate. This contradicts our assumptions about~$\vS$.
 %
\end{proof}

In view of Lemma~\ref{Ostars}, we shall call the subsets of a tree set~$\vS$ that split it the \emph{splitting stars} of~$\vS$.%
   \COMMENT{}

\section{Tangles in abstract separation systems}\label{sec:tangles}

In~\cite{TangleTreeAbstract, TangleTreeGraphsMatroids}, the tangles of graphs introduced by Robertson and Seymour~\cite{GMX} were generalized to arbitrary separation systems~$\vS$ as something called `$\F$-tangles'. Here, $\F$~is a collection of subsets of~$\vS$, usually stars, and an $\F${\em -tangle\/} is a consistent orientation of~$\vS$ that has no subset in~$\F$. 

When $\vS$ is a separation system of some combinatorial structure, such as a graph, such orientations are often defined as those that point to some `highly cohesive substructure', or~{\em HCS}, such as a large grid minor. These orientations will avoid any~$\F$ whose stars~$F$ each point to an area of the graph that is too small to accommodate such an~HCS, and thus will be $\F$-tangles.

The converse, that {\em all\/} $\F$-tangles point to an HCS of that kind, is often also true, at least qualitatively, but much harder to prove: those HCSs usually have, by definition, lots of concrete features that may have little to do with the fact that they are highly cohesive, and which are therefore difficult to reconstruct from only that assumption.%
   \COMMENT{}

The shift of paradigm that tangles have brought to the connectivity theory of graphs, now, is twofold. First, that when tangles are induced by concrete HCSs as above, we can often prove important theorems about the latter just in terms of those tangles: we do not need to know exactly `what' those substructures are, only `where' they are. The second is that $\F$-tangles can be interesting structural objects even when they are not induced by a concrete~HCS: they seem to capture the essence of high cohesion better than the various concrete HCSs that each come with their own frills. This transition is explained further in the introduction of~\cite{TangleTreeAbstract}.

Still, the elements of~$\F$ will typically be defined in terms of the structure, such as a graph, that is being separated; see~\cite{TangleTreeGraphsMatroids, MonaLisa} for examples. But what if $\vS$ comes as an abstract separation system? Is it possible to define a `generic' kind of~$\F$ such that the $\F$-tangles of a separation system~$\vS$ will {\em always\/} describe highly cohesive substructures, in any structure which $\vS$ might separate?

This is indeed possible. Given a universe~$\vU$ of separations, $\T = \T(\vU)$ be the set of all subsets of~$\vU$ that consist of up to three separations whose supremum in $\vec U$ is co-small, i.e., has a small inverse.%
   \COMMENT{}

\begin{DEF}
A {\em tangle} of a separation system~$\vS\sub\vU\!$ is a $\T$-tangle of~$\vS$.
\end{DEF}

In other words, a tangle of~$\vS$ is a consistent orientation of~$\vS$ that has no subset $\{\vr,\vs,\vt\}$ such that $\rv\land\sv\land\tv$ is small. Note that $\rv\land\sv\land\tv$ is not required to lie in~$\vS$, only in~$\vU$. One can show that if $\vU$ is distributive as a lattice, and $\vS$ is submodular (see below), then the consistent orientations of $S$ without a star in~$\cal T$ have no other subsets in~$\cal T$ either, i.e., are tangles of~$S$.

The fact that these separations $\rv\land\sv\land\tv$ are small bears out the intuition that the sets $F\in\F$ for an $\F$-tangle should point to somewhere `small' in any structure that $\vS$ separates: in all known applications, small separations do indeed point away from an area that is intuitively `small'. What is new here is that this has been expressed purely in terms of~$\vS$.

For graphs, these abstract tangles differ slightly from the original graph tangles of Robertson and Seymour~\cite{GMX}, but the two are closely related.%
   \COMMENT{}
   For matroids, they yield precisely the familiar matroid tangles~-- including the `principal' ones oriented towards some fixed elements of the matroid, which are sometimes excluded from consideration.%
   \COMMENT{}

For graph and matroid tangles, there are two fundamental types of theorem: {\em tangle-tree theorems\/} which assert the existence of a \td\ of the graph or matroid whose separations (those associated with the edges of the decomposition tree) distinguish all the tangles considered pairwise, and {\em tangle duality\/} theorems which say that any graph or matroid thas has no such tangle must have a \td\ that witnesses this, in the sense that its parts are too small to accommodate such a tangle.

These two theorems have been proved for abstract separation systems too: the tangle-tree theorem was proved in~\cite{ProfilesNew} for {\em profiles}~-- a~notion generalising tangles~-- while the tangle duality theorem was proved in~\cite{TangleTreeAbstract} for $\F$-tangles whose $\F$ satisfies some technical conditions.

The good news, now, is that for our generic `abstract tangles' we have both a tangle-tree theorem and a tangle duality theorem.

To state this concisely, we need a few more definitions for a separation system $\vS$ in a universe~$\vU$. Let us say that a tree set $\tau\sub\vS$ {\em distinguishes\/} a set $\O$ of orientations of~$\vS$ if for all distinct $O,O'\in\O$ there exists some $\vs\in\tau$ such that $O$ and~$O'$ orient $s$ differently. And let us call $\tau$ a tree set {\em over} a set $\F\sub 2^\vS\!$ of stars if all its splitting stars lie in~$\F$.
Finally, we call $\vS$ {\em submodular\/} in~$\vec U$ if for every two elements of~$\vec S$ their infimum or their supremum in~$\vec U$ also lies in~$\vec S$.

\begin{THM} {\rm\cite{AbstractTangles}}
Every submodular separation system $\vS$ in a universe~$\vU\!$ of separations contains a tree set that distinguishes all its tangles.
   \COMMENT{}
\end{THM}

There will likely also be a duality theorem for abstract tangles in~\cite{AbstractTangles}.%
   \COMMENT{}

\bibliographystyle{plain}
\bibliography{collective}

\begin{thebibliography}{10}

\bibitem{duality1inf}
N.~Bowler, R.~Diestel, and F.~Mazoit.
\newblock Tangle-tree duality in infinite graphs.
\newblock In preparation.

\bibitem{confing}
J.~Carmesin, R.~Diestel, F.~Hundertmark, and M.~Stein.
\newblock Connectivity and tree structure in finite graphs.
\newblock {\em Combinatorica}, 34(1):1--35, 2014.

\bibitem{TreeSets}
R.~Diestel.
\newblock Tree sets.
\newblock To appear in {\em Order} (2017), DOI 10.1007/s11083-017-9425-4,
  arXiv:1512.03781.

\bibitem{DiestelBook16}
R.~Diestel.
\newblock {\em Graph Theory \emph{(5th edition, 2016)}}.
\newblock Springer-Verlag, 2017.
\newblock \\ Electronic edition available at {\small\tt
  http://diestel-graph-theory.com/}.

\bibitem{AbstractTangles}
R.~Diestel and J.~Erde.
\newblock Tangles in abstract separation systems.
\newblock In preparation.

\bibitem{ProfileDuality}
R.~Diestel, J.~Erde, and Ph. Eberenz.
\newblock Duality theorem for blocks and tangles in graphs.
\newblock arXiv:1605.09139, to appear in {\em SIAM J.~Discrete Math}.

\bibitem{ProfilesNew}
R.~Diestel, F.~Hundertmark, and S.~Lemanczyk.
\newblock Profiles of separations: in graphs, matroids, and beyond.
\newblock arXiv:1110.6207, to appear in {\em Combinatorica}.

\bibitem{TangleTreeAbstract}
R.~Diestel and S.~Oum.
\newblock Tangle-tree duality in abstract separation systems.
\newblock arXiv:1701.02509, 2017.

\bibitem{TangleTreeGraphsMatroids}
R.~Diestel and S.~Oum.
\newblock Tangle-tree duality in graphs, matroids and beyond.
\newblock arXiv:1701.02651, 2017.

\bibitem{MonaLisa}
R.~Diestel and G.~Whittle.
\newblock Tangles and the {M}ona {L}isa.
\newblock arXiv:1603.06652.

\bibitem{GMX}
N.~Robertson and P.D. Seymour.
\newblock Graph minors. {X}. {O}bstructions to tree-decomposition.
\newblock {\em J.~Combin.\ Theory (Series B)}, 52:153--190, 1991.

\end{thebibliography}

\small\vfill\noindent Version 17.4.2017

\end{document}